\documentclass{amsart}
\usepackage{graphicx}
\usepackage{graphicx}
\usepackage[pdfauthor={Richard J. Mathar},pdfkeywords={Domino tilings, Polyominoes, Combinatorics},pdfsubject={Mathematics, Combinatorics},colorlinks=true,citecolor=blue]{hyperref}

\usepackage{amsmath}

\newtheorem{remark}{Remark}
\begin{document}

\title[Tiling Rectangular Regions with Rectangular Tiles]{Tilings of Rectangular Regions by Rectangular Tiles: Counts
Derived from Transfer Matrices.}

\author{Richard J. Mathar}
\email{mathar@mpia.de}
\urladdr{http://www.mpia.de/~mathar}
\address{Hoeschstr. 7, 52372 Kreuzau, Germany}

\subjclass[2010]{Primary 52C20, 05B45; Secondary 05A15}

\date{\today}
\keywords{Tiling, Rectangular Mats, Combinatorics}

\begin{abstract}
Step by step completion of a left-to-right tiling of a rectangular floor with tiles of a single
shape starts from one edge of the floor, considers  the possible ways of inserting
a tile at the leftmost uncovered square, passes through a sequence of rugged shapes of the front
line between covered and uncovered regions of the floor, and finishes with a straight front line
at the opposite edge. We count the tilings by mapping the front shapes to 
nodes in a digraph, then counting closed walks on that digraph with the transfer matrix method.

Generating functions are detailed for tiles of shape $1\times 3$, $1\times 4$ and $2\times 3$ and modestly
wide floors. Equivalent results are shown for the 3-dimensional analog of filling bricks
of shape $1\times 1\times 2$, $1\times 1\times 3$, $1\times 1\times 4$, $1\times 2\times 2$
or $1\times 2\times 3$
into rectangular containers of small cross sections.
\end{abstract}

\maketitle 

\section{Definitions (Dimensions, State Vectors)} 
Given a floor of width $m$ and length $n$ and a prototile of width $t_m$ and length $t_n$, $t_m\le t_n$,
we consider the number of ways of covering the floor by
\begin{equation}
N=\frac{mn}{t_nt_m}
\label{eq.areaDiv}
\end{equation}
non-overlapping tiles.
Tiles may be placed on the floor mixing both orientations. Only the case of
coprime $t_n$ and $t_m$ is of interest, because the geometry could
otherwise be shrunk by the common factor without changing the number of tilings.

The symmetry of the $m\times n$ rectangle will not be taken into account;
tilings which are equivalent to other tilings through reflections or rotations of
the stack are counted including their multiplicity due to their (missing) symmetry.

\begin{figure}
\caption{
Illustration of the vector of heights, the profile of the covered region
into the uncovered part of the floor.
}
\includegraphics[width=0.6\textwidth]{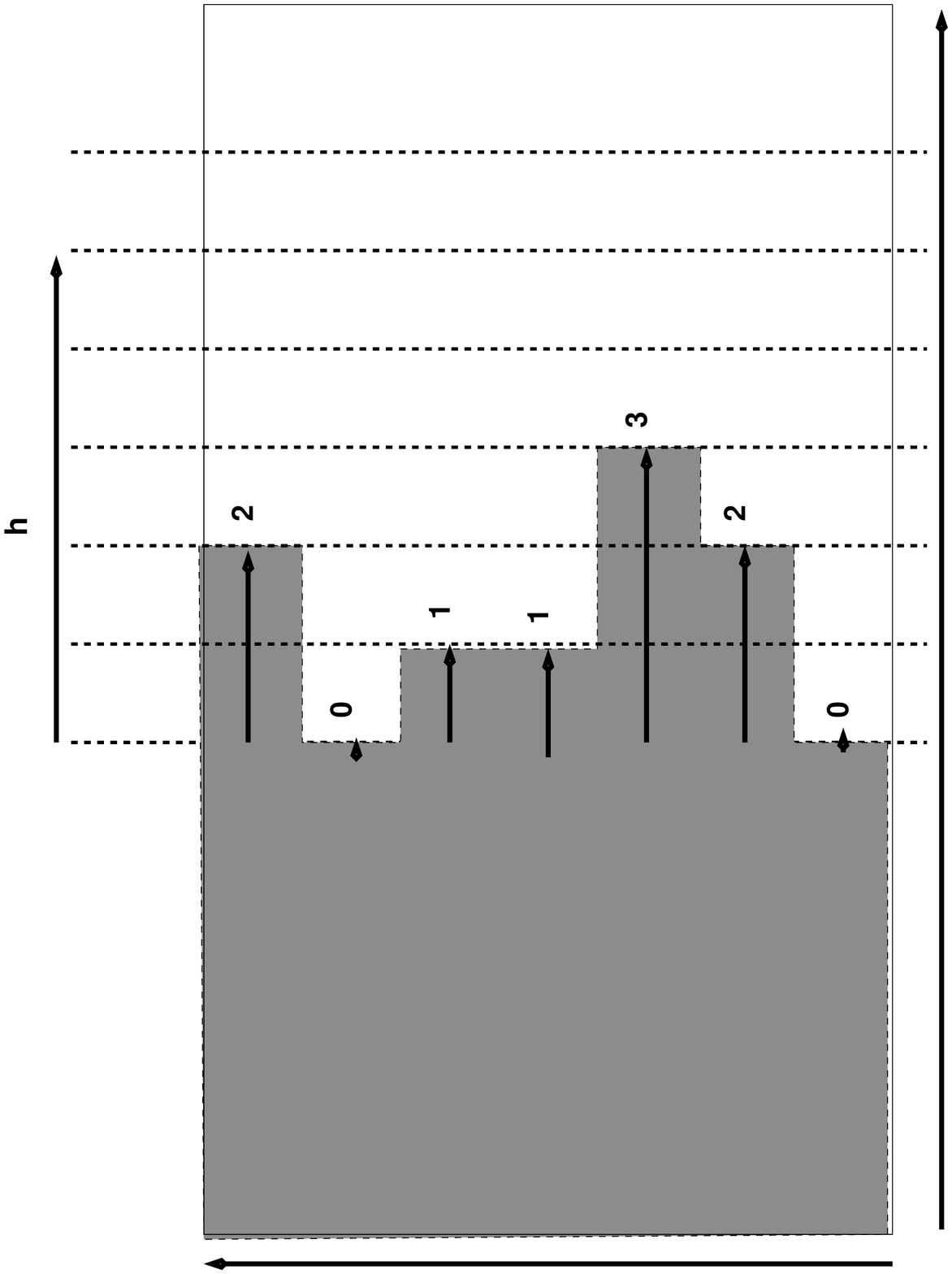}
\label{fig.heig}
\end{figure}

The following counting technique reflects the act of paving
the floor starting from the (left) short edge.
The rule
of placing the next tile is to cover the leftmost not yet covered unit square
of the floor, and if there are more than one of them to cover, the one closest
to the front edge of the floor.
Each intermediate tiling shall
be  characterized by a \emph{state vector} $(h_1h_2\ldots h_m)$ of heights $0\le h_i\le t_n$
that encode how far the tiles at vertical position $1\le i\le m$ penetrate into the
uncovered territory beyond the rightmost occupied square \cite{ReadFQ18,ReadAM24}.
Figure \ref{fig.heig} illustrates these integers $h_i$
of a state vector $(0231102)$
with a floor width of $m=7$.

A directed graph is defined where the different state vectors are the
nodes; edges point from each state vector
to the state vectors that can be reached by placing the next tile,
so walking along a single edge represents adding a tile.
The possibility of reaching a state where no further tile can be placed
in conjunction with the two orientations of the tile
lead to outdegrees in the range 0 up to 2.

The number of state vectors (of nodes in the graph) is finite.
An upper limit is $(1+t_n)^m$, because the rule of the
next placement puts a limit of $t_n$, the tile's long edge, to each profile.
In addition, the maximum height $h_i=t_n$ may only appear for the first consecutive
vector elements starting at $i=1$, because otherwise at least one tile would have been placed
prematurely at an advanced position reserved for later placements.
Further reduction of the upper limit stems from the requirement that
at least one $h_i$ must equal zero.

\section{Counting Closed Walks on the Digraph}
\subsection{Generating Function}
We wish to count the number $T(N)$ of different closed walks of $N$ steps starting from
and returning to the state vector $(00\ldots00)$, which represents the straight empty front
edge at the start of the paving and represents also the straight front edge when the opposite
edge is reached. The aforementioned rule to cover the ``lowest'' unoccupied unit square
(employing some sequential enumeration of the $nm$ unit squares of the floor) ensures
that there is a 1-to-1 map of the different tilings to the walks on the graph.

The number of walks are registered by construction of their ordinary generating function
(GF)
\begin{equation}
T_{t_m\times t_n}(m,z)\equiv \sum_{N\ge 0} T_{t_m\times t_n}(m,N) z^N.
\end{equation}
\begin{remark}
This definition uses the tile count $N$ to manage the results, whereas my earlier
notation uses the floor length $n$ to tabulate them \cite{MatharArxiv1311}. For the cases studied here,
the shape $t_m\times t_n$ of the tile and the width $m$ of the floor are kept fixed,
so $N$ and $n$ are easily translated using \eqref{eq.areaDiv}.
\end{remark}
Some relation between the GF and the
graph topology can easily be established akin to the Kirchhoff rules of impedances
of circuit networks \cite{MurataIEEE77}:
\begin{itemize}
\item
(Disconnected) Subgraphs not connected to $(00\ldots0)$ are discarded.
Their contribution to the generating
function is multiplication by the factor 1.
\item
(Edge weights) A step along a single edge multiplies the generating function by $z$.
\item
(Serial Paths) A non-splitting chain of walks multiplies the generating functions along
the sub-paths.
\item
(Parallel Paths) At a node with outdegree larger than one, where diverting alternative walks
lead across disconnected subgraphs and rejoin later, so
exactly one of the paths is to be taken to reach the join, the sum
of the generating functions of the dispatched subgraphs is build.
\end{itemize}
These basics may be combined to formulate rules for
simply weaved networks:
\begin{itemize}
\item
A circuit from a node back to itself, which allows to walk
that circuit any number of times or not at all, contributes $1/(1-l(z))$ to
the generating function, where $l(z)$ is the GF of the circuit.
The formula sums a geometric series
of the two rules for chains and parallel paths.
\item
A node with several disconnected walks returning to the node, where the
generating functions are $T^{(1)}(z)$, $T^{(2)}(z)$ etc.\ and where walks may be executed
in any order, represents
\begin{equation}
\frac{1}{1-\sum_i T^{(i)}(z)}=1+\sum_i T^{(i)}(z)
+ \sum_{ij} T^{(i)}(z)T^{(j)}(z)+\cdots.
\label{eq.parLoop}
\end{equation}
The terms on the right hand side accumulate (as products) the number of walks
choosing first circuit $i$, then circuit $j$ etc, each as many times as wished.
\end{itemize}

\begin{figure}
\caption{
State diagram while tiling $2\times n$ floors with
dominoes, generating the Fibonacci sequence.
}
\includegraphics[width=0.18\textwidth]{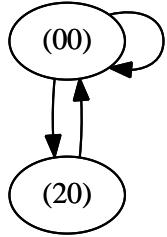}
\label{fig.state221}
\end{figure}

Figure \ref{fig.state221} shows the fundamental example of placing
$1\times 2$ tiles on $2\times n$ floors. Orienting the tile parallel
to $n$ leads from $(00)$ to $(20)$ or (in the second lane) from $(20)$
to $(00)$, whereas orienting it parallel to the short edge loops from $(00)$
back to $(00)$.
The disconnected subgraph with the circuit $(01)\to (10)\to (01)$ is not plotted.
Starting at $(00)$ there is one loop with a single step, with GF $z$.
The circuit passing by the intermediate $(20)$ represents a walk with two steps,
with GF $z^2$.
The rule \eqref{eq.parLoop} for multiple circuits attached to a
node combines these two GF's to $T_{1\times 2}(2,z)=1/(1-z-z^2)$, the Fibonacci series
\cite[A000045]{EIS}.

\begin{figure}
\caption{
State diagram while tiling $3\times n$ floors with
dominoes.
}
\includegraphics[width=0.5\textwidth]{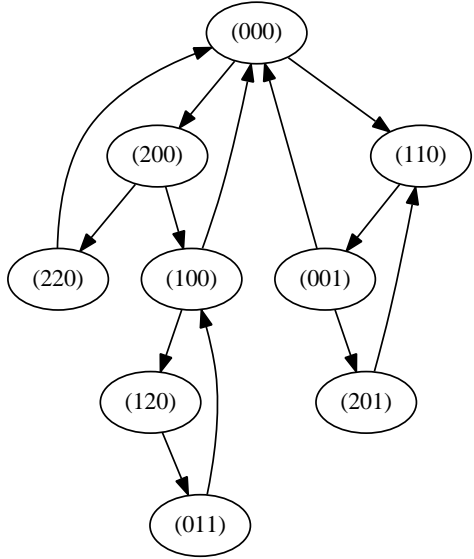}
\label{fig.state321}
\end{figure}

A second example of application is the tiling of $3\times n$ floors
with $1\times 2$ tiles, Figure \ref{fig.state321}.
One option is to walk $(000)\rightarrow (200)\rightarrow (220)\rightarrow (000)$
which contributes $T^{(1)}=z^3$ by the chain rule.
There is a subwalk $(100)\rightarrow (120)\rightarrow (011)
\rightarrow (100)$ in a circuit which represents $1/(1-z^3)$; since it is reached from
$(000)$ via $(200)$ with two steps and returns to $(000)$
in one step, it actually represents $T^{(2)}=z^3/(1-z^3)$. Finally there is 
a walk to $(110)\rightarrow (001)$ with an option to circle any number of times
via $(201)\rightarrow (110)\rightarrow (001)$ before moving back to $(000)$,
$T^{(3)}=[z^2/(1-z^3)]z$. The effective total of the three circuits
originating from $(000)$ is
\cite[A001835]{EIS}
\begin{equation}
T_{1\times 2}(3,z)=\frac{1}{1-T^{(1)}-T^{(2)}-T^{(3)}}
=\frac{(1-z)(1+z+z^2)}{1-4z^3+z^6}
=1+3z^3+11z^6+41z^9+\cdots.
\label{eq.state321}
\end{equation}

\subsection{Irreducible GF}
The function $T(m,N)$ counts tilings including any number of blocks
of side-by-side smaller tilings $N=N'+N''+\cdots$ of the common width $m$.
In the speak of the digraph, $T(m,z)$ counts walks that pass through $(00\ldots 0)$
any number of times. The \emph{irreducible} tilings may be defined as those
that cannot be cut into smaller tilings by cuts of length $m$ parallel to the
short side of the floor. Their walk returns to $(00\ldots 0)$ just once.
Their GF shall be denoted by $\hat T(m,z)$:
\begin{equation}
T(m,z)=\frac{1}{1-\hat T(m,z)}.
\end{equation}
The reduction of \eqref{eq.state321} to these irreducible cases yields
for example
\begin{equation}
\hat T_{1\times 2}(3,z)=\frac{z^3(3-z^3)}{1-z^3}
=3z^3+2z^6+2z^9+2z^{15}+2z^{18}+\cdots.
\end{equation}

\subsection{Inverting Transfer Matrices}
An iterative procedure to compute the Taylor Series of the generating function
according to the rules of the previous section is to ``load'' the node $(00\ldots 0)$
with $T=1$ and all other nodes with $T=0$, and to multiply the vector of loads
with increasingly higher powers of the Transfer Matrix of the digraph. The Transfer
Matrix is defined to contain the weight $z$ of the edge at the column
associated with the start of the step and the row associated with the end of the step along
the edge, i.e., essentially the incidence matrix pairing direct predecessors and successors in the graph.
The number of $z$ in the columns equals the outdegree of the node;
the number of $z$ in the columns equals the indegree of the node.

The $9\times 9$ Transfer Matrix of Figure \ref{fig.state321} with 9 nodes
and a vector of 9 states sorted as
$(000)$,
$(001)$,
$(011)$,
$(100)$,
$(110)$,
$(120)$,
$(200)$,
$(201)$,
and $(220)$
becomes
\begin{equation}
X=
\left(
\begin{array}{rrrrrrrrr}
0 & z & 0 & z & 0 & 0 & 0 & 0 & z \\ 
0 & 0 & 0 & 0 & z & 0 & 0 & 0 & 0 \\ 
0 & 0 & 0 & 0 & 0 & z & 0 & 0 & 0 \\ 
0 & 0 & z & 0 & 0 & 0 & z & 0 & 0 \\ 
z & 0 & 0 & 0 & 0 & 0 & 0 & z & 0 \\ 
0 & 0 & 0 & z & 0 & 0 & 0 & 0 & 0 \\ 
z & 0 & 0 & 0 & 0 & 0 & 0 & 0 & 0 \\ 
0 & z & 0 & 0 & 0 & 0 & 0 & 0 & 0 \\ 
0 & 0 & 0 & 0 & 0 & 0 & z & 0 & 0 \\ 
\end{array}
\right)
\end{equation}
for example. The two $z$ in its first column instantiate the walks $(000)\to(110)$
and $(000)\to(200)$.

The full GF follows by (i) summation of the 
geometric series of the matrix powers, which ultimatively is the inversion
$T(z)=\sum_{i\ge 0}[X(z)]^i = \mathrm{inv} (I-X(z))$ with $I$ the identity matrix,
(ii) multiplying that inverse
from the right by the initial load $(1,0,0,\ldots)$ and (iii) taking the element of that
vector associated with $(00\ldots 0)$ \cite{Stanley}.

\begin{remark}
Our application does not rely on the full inverse of the Transfer Matrix but only
its first column; the effort reduces to solving a linear system of equations.
As the number of states is finite, solving the equation by
Cramer's rule proves that the GF's are rational functions of $z$.
\end{remark}

\section{Floor Tilings}

\subsection{Dominoes}

Tiling with dominoes is much better studied in the literature
than tiling with other polyominoes because thermodynamic properties
of dimer coverings are of interest to the physics of surfaces
\cite{KlarnerDM32,HockDAM8,StanleyDAM12,StrehlAAM27,RuskeyEJC16,LiebJMP8}.
We skip this special subject because no new results arise in the
present context. For increasing width $m$, the sequences $T_{1\times 2}(m,n)$
are entries A000075, A005178, A003775, A028469--A028474 in
the Encyclopedia of Integer Sequences \cite{EIS}.

\subsection{$1\times 3$ Straight Trominoes}
\begin{figure}
\caption{
State diagram while tiling $3\times n$ floors with
$1\times 3$ 3-ominoes.
}
\includegraphics[width=0.22\textwidth]{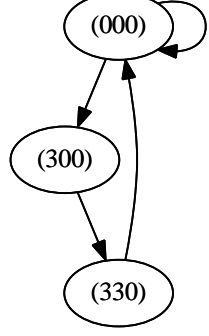}
\label{fig.state331}
\end{figure}

For tiles of $1\times 3$ shape we find
\begin{equation}
T_{1\times 3}(2,z)= \frac{1}{1-z^2},
\label{eq.state231}
\end{equation}
which essentially says that there is one tiling whenever \eqref{eq.areaDiv} is an integer,
because the long side of the tile needs to be aligned with the long side of the floor.
We find via Figure \ref{fig.state331}
\cite[A000930]{EIS}
\begin{equation}
T_{1\times 3}(3,z)=
\frac{1}{1-z-z^3}
=
1+z+z^{2}+2z^{3}+3z^{4}+4z^{5}+6z^{6}+9z^{7}+13
z^{8}+19z^{9}+28z^{10}+\cdots.
\label{eq.state331}
\end{equation}

\begin{figure}
\caption{
State diagram while tiling $4\times n$ floors with
$1\times 3$ 3-ominoes.
}
\includegraphics[width=0.7\textwidth]{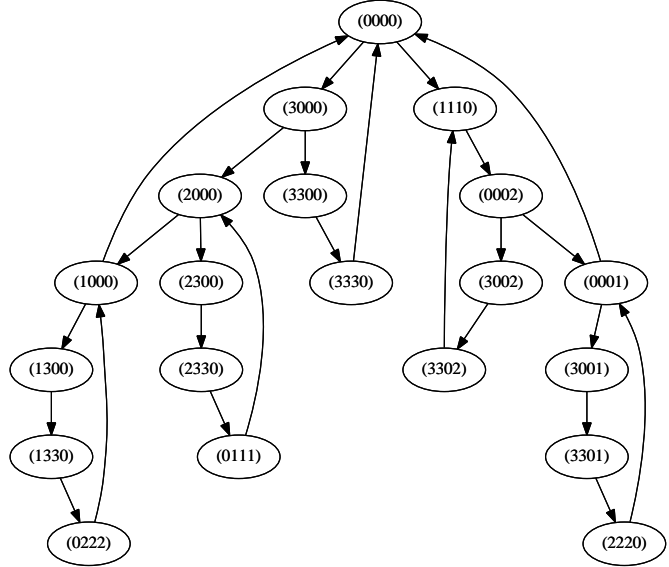}
\label{fig.state431}
\end{figure}

For floors of width $m=4$ the digraph of Figure \ref{fig.state431}
reduces to sequence
\cite[A049086]{EIS}
and basically the even numbers:
\begin{equation}
T_{1\times 3}(4,z)=
{\frac { \left( -1+z \right) ^{2} \left( z+1 \right) ^{2} \left( z^{
2}+1 \right) ^{2}}{-5z^{4}+1+3z^{8}-z^{12}}}
=
1+3z^{4}+13z^{8}+57z^{12}+249z^{16}+1087z^{20}+\cdots;
\end{equation}
\begin{equation}
\hat T_{1\times 3}(4,z)=
\frac {z^{4} \left( 3-2z^{4}+z^{8} \right) }
{ 1-2z^4+z^8}
=
3z^{4}+4z^{8}+6z^{12}+8z^{16}+10z^{20}+12z^{24
}+14z^{28}+16z^{32}+\cdots.
\end{equation}

\begin{figure}
\caption{
State diagram while tiling $5\times n$ floors with
$1\times 3$ 3-ominoes.
}
\includegraphics[width=0.8\textwidth]{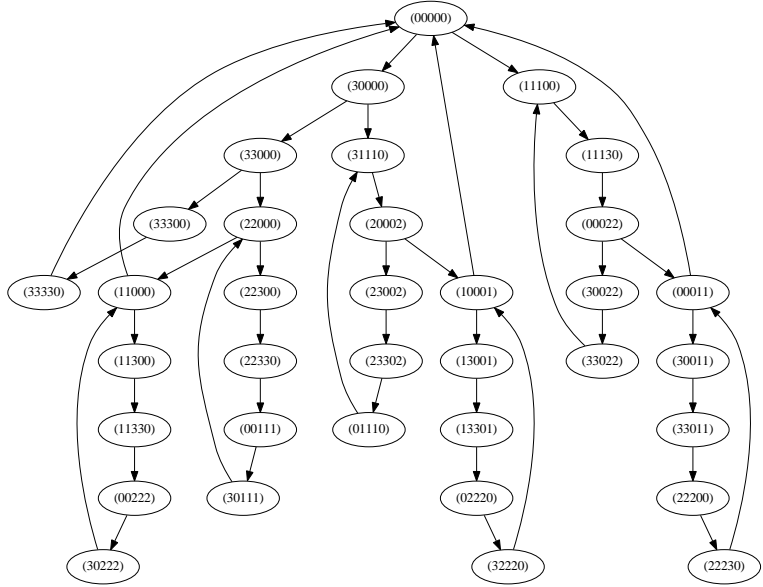}
\label{fig.state531}
\end{figure}

Figure \ref{fig.state531} gives \cite[A236576]{EIS}
\begin{equation}
T_{1\times 3}(5,z)=\frac{(1-z^5)^2}{1-6z^5+3z^{10}-z^{15}}
=
1+4z^{5}+22z^{10}+121z^{15}+664z^{20}+3643z^{25}+
19987z^{30}+\cdots
\end{equation}
and basically multiples of three:
\begin{equation}
\hat T_{1\times 3}(5,z)=
\frac {z^{5} \left( 4-2z^{5}+z^{10} \right) }
{ (1-z^5)^2}
=
4z^{5}+6z^{10}+9z^{15}+12z^{20}+15z^{25}+18z^{
30}+21z^{35}+24z^{40}+27z^{45}+\cdots.
\end{equation}

\begin{figure}
\caption{
State diagram while tiling $6\times n$ floors with
$1\times 3$ 3-ominoes.
}
\includegraphics[width=0.7\textwidth]{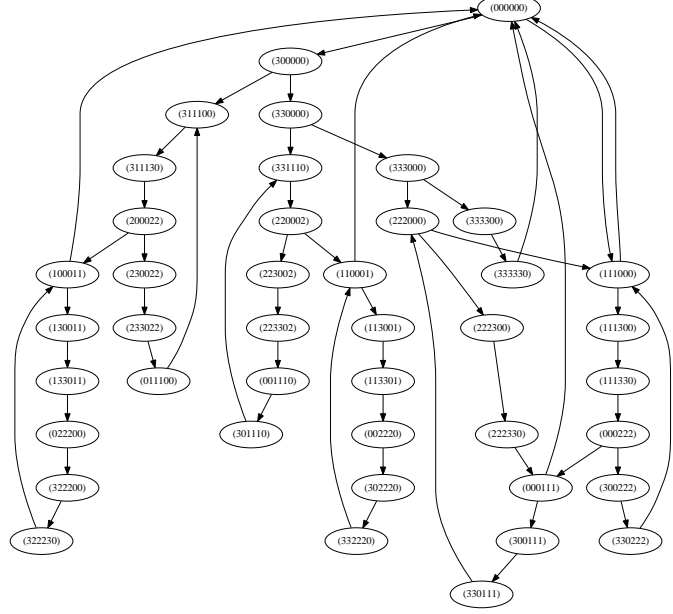}
\label{fig.state631}
\end{figure}

Figure \ref{fig.state631} reduces to \cite[A236577]{EIS}
\begin{multline}
T_{1\times 3}(6,z)=
{\frac { \left( 1-z^{6} \right) ^{2} \left(-z^{4}+1-z^{6}
 \right) }{-z^{20}+z^{24}+z^{22}+10z^{12}-5z^{18}-3z^
{16}+z^{14}+z^{8}-7z^{6}+5z^{10}-z^{4}-z^{2}+1}}
\\
=
1+z^{2}+z^{4}+6z^{6}+13z^{8}+22z^{10}+64z^{12}+155
z^{14}+321z^{16}+783z^{18}+1888z^{20}+4233z^{22}+\cdots.
\end{multline}

The Transfer Matrix for filling $7\times n $ floors with $1\times 3$ 3-ominoes
has dimension $273\times 273$ \cite[A236578]{EIS},
\begin{equation}
T_{1\times 3}(7,z)=\frac{p_{1\times 3}(7,z)}{q_{1\times 3}(7,z)}
=
1+9z^{7}+155z^{14}+2861z^{21}+52817z^{28}+972557z
^{35}+\cdots
\end{equation}
with numerator
\begin{multline}
p_{1\times 3}(7,z)\equiv \left( z^{7}-1 \right) ^{2} (-z^{105}+14z^{98}-
104z^{91}+527z^{84}-1971z^{77}+5573z^{70}-11973z^{63}
\\
+19465z^{56}-23695z^{49}+21166z^{42}-13512z^{35}+
5915z^{28}-1685z^{21}+291z^{14}-27z^{7}+1);
\end{multline}
and denominator
\begin{multline}
q_{1\times 3}(7,z)\equiv 
-17
z^{119}+293180z^{56}-236178z^{49}+142400z^{42}-62621{
z}^{35}+19420z^{28}
\\
-4062z^{21}+533z^{14}-38z^{7}+z^{
126}+1+151z^{112}-946z^{105}+4558z^{98}-17135z^{91}
\\
+50164z^{84}
-114198z^{77}+202080z^{70}-277277z^{63}.
\end{multline}
\begin{equation}
\hat T_{1\times 3}(7,z)
=
9z^{7}+74z^{14}+800z^{21}+8398z^{28}+85908z^{35}+
867148z^{42}+8697028z^{49}+86962830z^{56}+\cdots.
\end{equation}

Filling $8\times n $ floors with $1\times 3$ 3-ominoes
is counted by
\begin{equation}
T_{1\times 3}(8,z)=\frac{p_{1\times 3}(8,z)}{q_{1\times 3}(8,z)}
=
1+13z^{8}+321z^{16}+8133z^{24}+204975z^{32}+5158223
z^{40}\cdots,
\end{equation}
with
\begin{multline}
p_{1\times 3}(8,z)\equiv
133678z^{32}+1-190305075z^{88}+2914z^{192}-19827 z^{184}
-326z^{200}+223054092z^{96}
\\
+108161z^{176}-486843 z^{168}
+25z^{208}+35967130z^{64}-51z^{8}-z^{216}+ 3806952z^{48}
-13264117z^{56}
\\
-15217z^{24}-5879363z^{152} +179445867z^{112}
-218171819z^{104}+1838302z^{160}+16010861 z^{144}
\\
+134225178z^{80}
-77378505z^{72}+73831771z^{128}-
124861824z^{120}+1148z^{16}-37212426z^{136}-830622z^{40
};
\end{multline}
\begin{multline}
q_{1\times 3}(8,z)\equiv 
233525z^{32}+1+z^{224}-609762885z^{88}+24393z^{192}-
146966z^{184}-3312z^{200}
\\
+779625485z^{96}+738848z^{176}
-3126151z^{168}+351z^{208}+88153581z^{64}-64z^{8}-26{
z}^{216}
\\
+7843386z^{48}-29769135z^{56}-24373z^{24}-33883064
z^{152}+739898086z^{112}-829644305z^{104}
\\
+11180105z^{
160}+87159919z^{144}+393536359z^{80}-207405453z^{72}+
353345037z^{128}
\\
-555959103z^{120}+1659z^{16}-190440779z
^{136}-1575184z^{40}.
\end{multline}

\subsection{$1\times 4$ Straight 4-ominoes}
Similar to \eqref{eq.state231}--\eqref{eq.state331},
widths $m$  smaller than the long edge of the polyomino lead to simple results
\cite[A003269]{EIS}:
\begin{equation}
T_{1\times 4}(3,z)=\frac{1}{1-z^3};
\end{equation}
\begin{equation}
T_{1\times 4}(4,z)=\frac{1}{1-z-z^4}
=1+z+z^{2}+z^{3}+2z^{4}+3z^{5}+4z^{6}+5z^{7}+7z
^{8}+10z^{9}+14z^{10}+19z^{11}
+\cdots .
\label{eq.state314}
\end{equation}

\begin{figure}
\caption{
State diagram while tiling $5\times n$ floors with
$1\times 4$ 4-ominoes.
}
\includegraphics[width=0.7\textwidth]{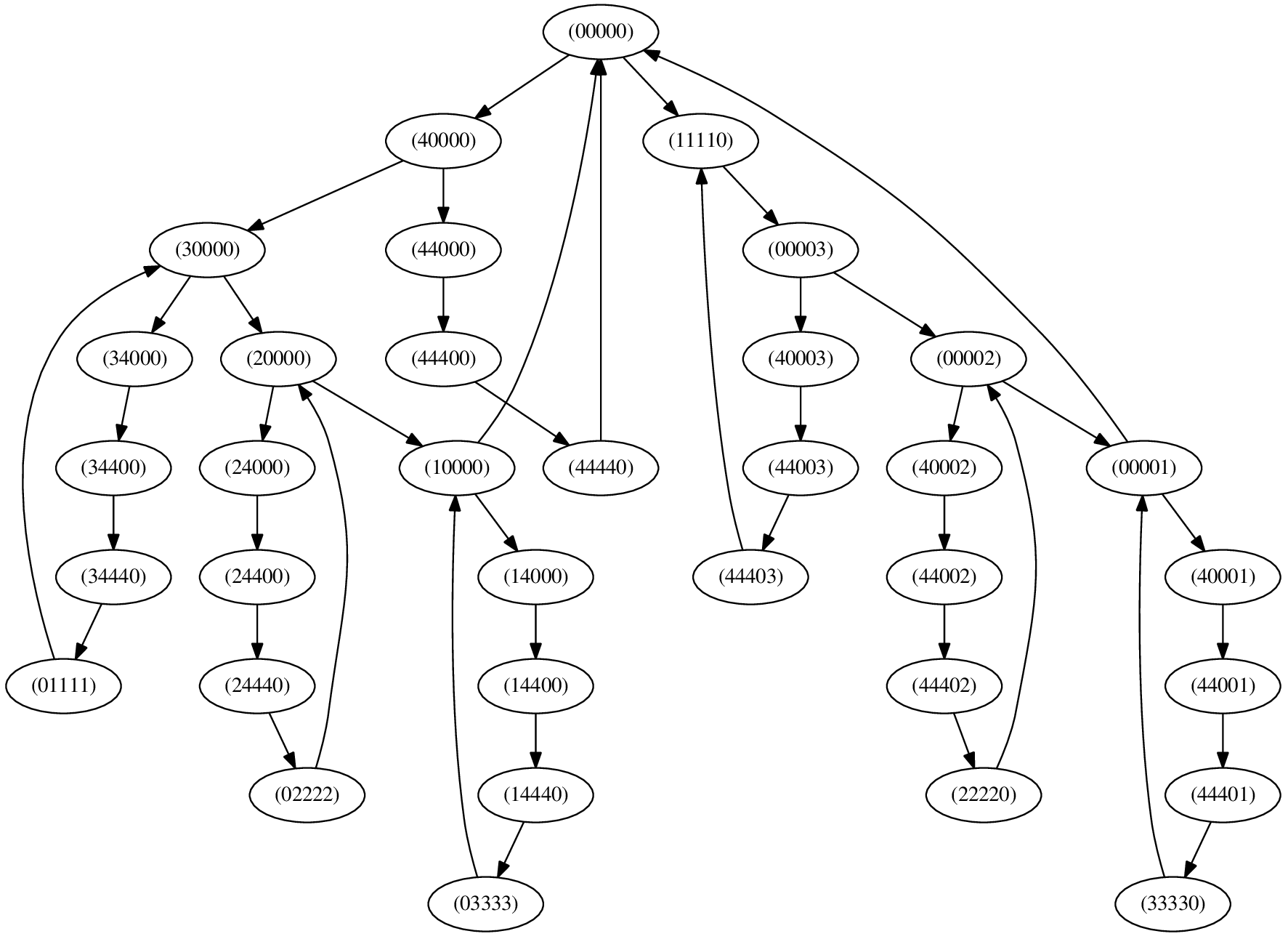}
\label{fig.state541}
\end{figure}

Figure \ref{fig.state541} is condensed to \cite[A236579]{EIS}
\begin{equation}
T_{1\times 4}(5,z)=
{\frac { \left( 1-z^{5} \right) ^{3}}{-6z^{5}+1+6z^{10}-4
z^{15}+z^{20}}}
=
1+3z^{5}+15z^{10}+75z^{15}+371z^{20}+1833z^{25}+\cdots
\end{equation}
and
\cite[A002378]{EIS}
\begin{equation}
\hat T_{1\times 4}(5,z)=
\frac {z^{5} \left( 3-3z^{5}+3z^{10}-z^{15} \right) }
{ (1-z^5)^3}
=
3z^{5}+6z^{10}+12z^{15}+20z^{20}+30z^{25}+42z^
{30}+56z^{35}+\cdots.
\end{equation}

\begin{figure}
\caption{
State diagram while tiling $6\times n$ floors with
$1\times 4$ 4-ominoes.
}
\includegraphics[width=0.9\textwidth]{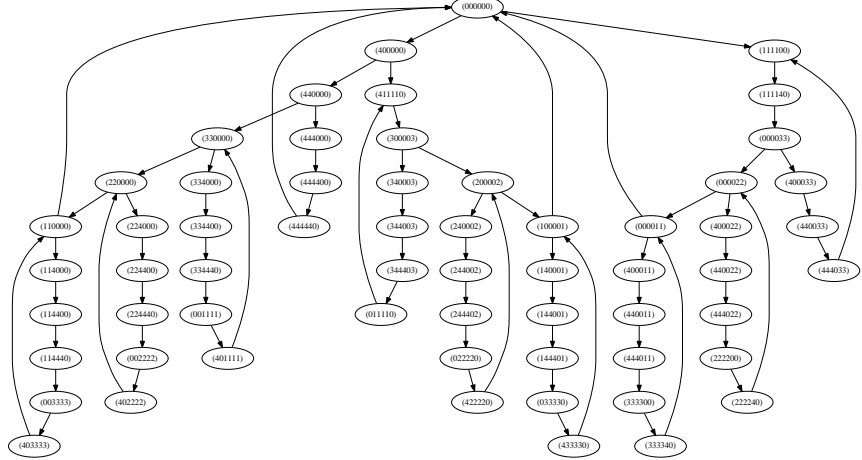}
\label{fig.state641}
\end{figure}

Figure \ref{fig.state641} represents \cite[A236580]{EIS}
\begin{equation}
T_{1\times 4}(6,z)=
{\frac { \left( 1-z^{6} \right) ^{3}}{-7z^{6}+1+6z^{12}-4
z^{18}+z^{24}}}
=
1+4z^{6}+25z^{12}+154z^{18}+943z^{24}+5773z^{30}+\cdots.
\end{equation}
Finally \cite[A236581]{EIS}
\begin{equation}
T_{1\times 4}(7,z)=
{\frac { \left( 1-z^7 \right) ^{3}}{-8z^{7}+1+6z^{14}-4z^{21}+z^
{28}}}
=
1+5z^{7}+37z^{14}+269z^{21}+1949z^{28}+14121z^{35}+\cdots.
\end{equation}
and \cite[A236582]{EIS}
\begin{equation}
T_{1\times 4}(8,z)=\frac{
p_{1\times 4}(8,z)}
{q_{1\times 4}(8,z)}
=
1+{z}^{2}+{z}^{4}+{z}^{6}+7\,{z}^{8}+15\,{z}^{10}+25\,{z}^{12}+37\,{z}
^{14}+100\,{z}^{16}+229\,{z}^{18}+\cdots,
\end{equation}
where
\begin{equation}
p_{1\times 4}(8,z)=
\left(1- z^2 \right) ^{3}\left( {z}^{2
}+1 \right) ^{3} \left( {z}^{4}+1 \right) ^{3} \left( {z}^{12}-{z}^{8}-
{z}^{6}-{z}^{4}+1 \right);
\end{equation}
\begin{multline}
q_{1\times 4}(8,z)=
-z^{4}-13z^{20}-5z^{36}+8z^{
12}-z^{2}-z^{40}-9z^{8}+16z^{16}-13z^{24}-2z^{38}+1
\\
+10z^{28}+5z^{14}+6z^{30}-6z^{22}+z^{44}+6z^{32}+
z^{34}+2z^{10}-2z^{26}.
\end{multline}

\subsection{$2\times 3$ Hexominoes}
\begin{figure}
\caption{
State diagram while tiling $5\times n$ floors with
$2\times 3$ 6-ominoes.
}
\includegraphics[width=0.6\textwidth]{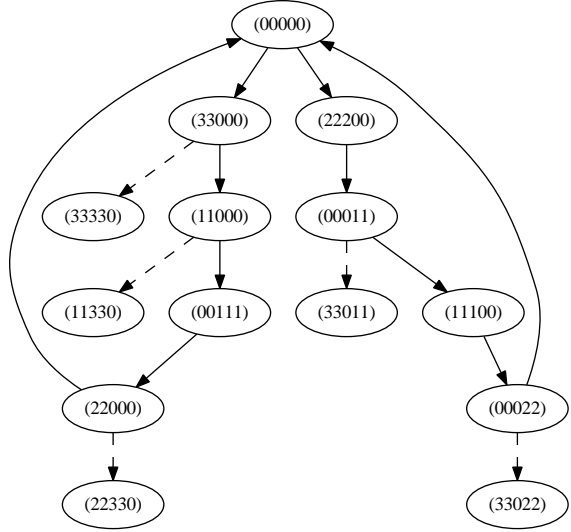}
\label{fig.state532}
\end{figure}

The state diagram in Figure \ref{fig.state532} is a first example with
``dangling'' nodes with zero outdegree. The shortest side $t_m=2$ of the
tile is larger than the slit of unit width if the height vector has the
pattern $(\ldots a0b\ldots)$ for some $a,b\ge 1$. After
pruning these dead ends (reducing the graph to the strongly connected subgraph),
the graph reduces to two isolated circuits of 5 edges
attached to the start position:
\begin{equation}
T_{2\times 3}(5,z)=
\frac{1}{1-2z^5}
=
1+2z^{5}+4z^{10}+8z^{15}+16z^{20}+32z^{25}+\cdots.
\end{equation}

\begin{figure}
\caption{
State diagram while tiling $6\times n$ floors with
$2\times 3$ 6-ominoes.
}
\includegraphics[width=0.4\textwidth]{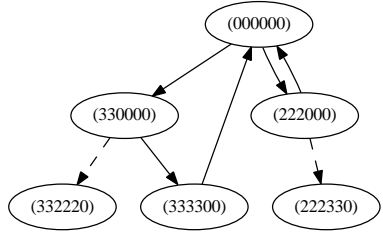}
\label{fig.state632}
\end{figure}

The state diagram in Figure \ref{fig.state632}
contains two nodes with zero outdegree, $(332220)$ and $(222330)$. After pruning
these, the topology is a circuit of 3 steps and a circuit of 2 steps
\cite[A182097]{EIS}:
\begin{equation}
T_{2\times 3}(6,z)=
\frac{1}{1-z^{2}-z^{3}}
=
1+z^{2}+z^{3}+z^{4}+2z^{5}+2z^{6}+3z^{7}+4z^{8}+
5z^{9}+7z^{10}+9z^{11}+\cdots.
\end{equation}

\begin{figure}
\caption{
State diagram while tiling $7\times n$ floors with
$2\times 3$ 6-ominoes.
}
\includegraphics[width=0.9\textwidth]{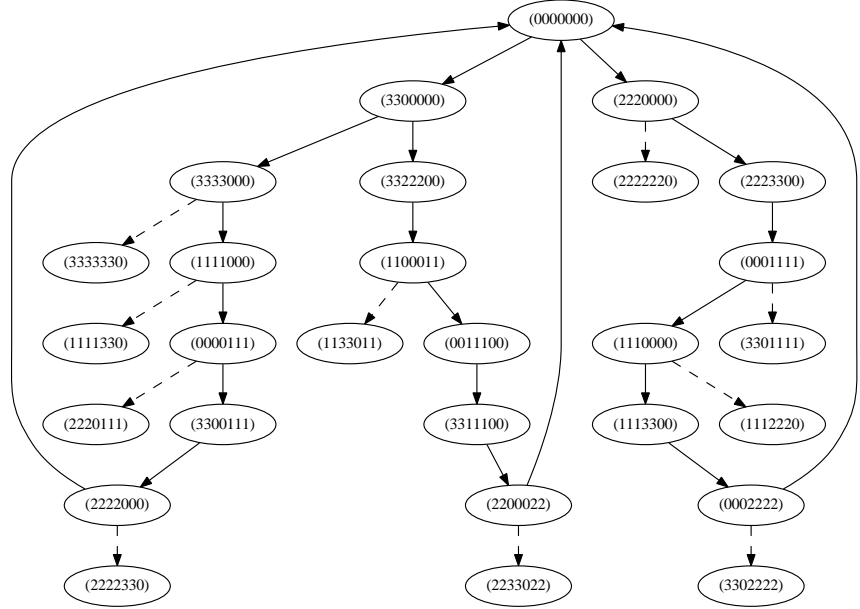}
\label{fig.state732}
\end{figure}

After pruning the dead ends of Figure \ref{fig.state732}, 
we are left with a skeleton of 3 circuits of 7 steps each, so
\cite[A000244]{EIS}
\begin{equation}
T_{2\times 3}(7,z)=
\frac{1}{1-3z^7} = 1+3z^6+9z^{14}+\cdots.
\end{equation}

\begin{figure}
\caption{
State diagram while tiling $8\times n$ floors with
$2\times 3$ 6-ominoes.
}
\includegraphics[width=0.99\textwidth]{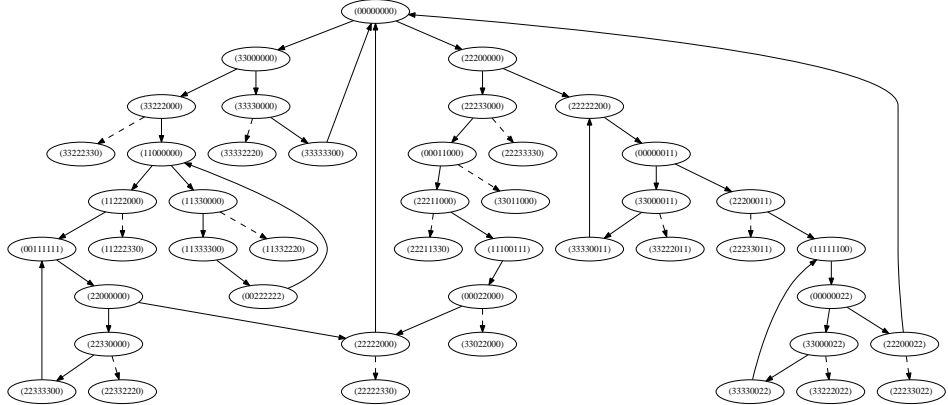}
\label{fig.state832}
\end{figure}

\begin{figure}
\caption{
State diagram while tiling $9\times n$ floors with
$2\times 3$ 6-ominoes.
}
\includegraphics[width=0.99\textwidth]{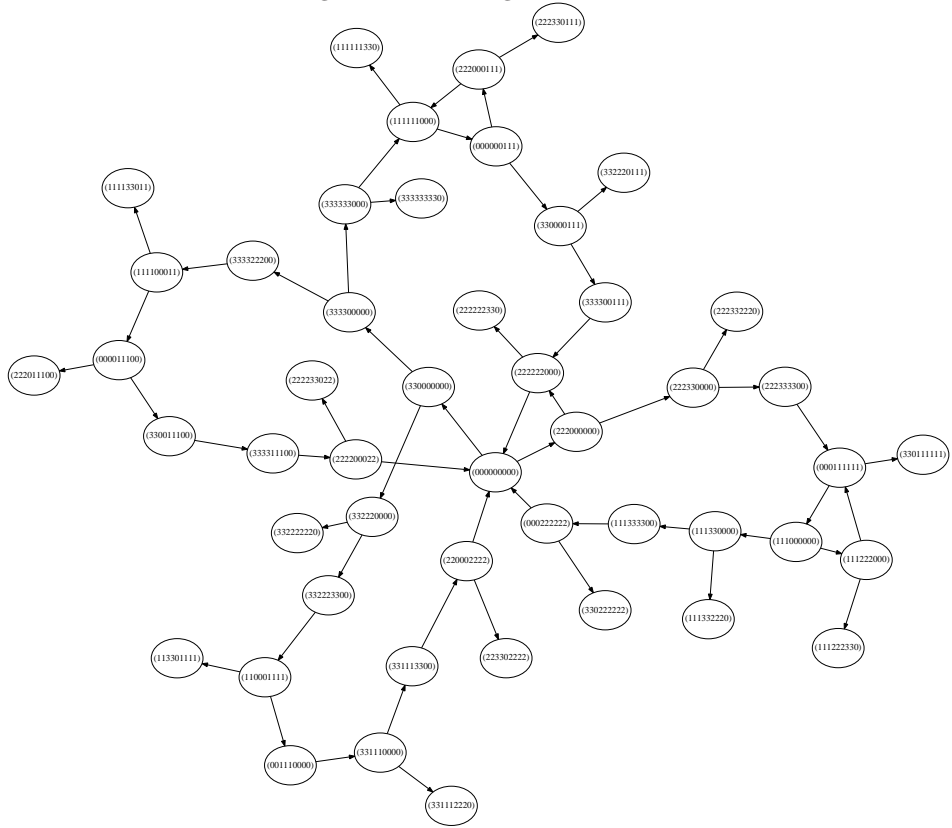}
\label{fig.state932}
\end{figure}

Figure \ref{fig.state832} reduces to \cite[A236583]{EIS}
\begin{equation}
T_{2\times 3}(8,z)=
{\frac { \left( -1+z^{4} \right) ^{2}}{z^{12}-z^{16}+1-3z^{
4}}}
= 
1+z^{4}+4z^{8}+11z^{12}+33z^{16}+96z^{20}+281z^{
24}+821z^{28}+\cdots,
\end{equation}
and the case of width $m=9$
\ref{fig.state932}
to \cite[A236584]{EIS}
\begin{equation}
T_{2\times 3}(9,z)=
{\frac {1-z^{3}}{-4z^{9}+1-2z^{3}+z^{6}+2z^{12}}}
=
1+z^{3}+z^{6}+5z^{9}+11z^{12}+19z^{15}+45z^{18}+
105z^{21}+\cdots.
\end{equation}
\clearpage

\section{Stacking Bricks in Rooms---The 3D analog}
\subsection{Digraphs}
The three-dimensional equivalent to the tilings packs
impenetrable bricks of dimension $t_k\times t_m\times t_n$ into containers
of shape $k\times m\times n$,  leaving no hole \cite{MincLMA8}. The volume
ratio $kmn/(t_kt_mt_n)$ is an integer.

No generic new aspects relative to the handling of the two dimensions
arise. Because the side lengths of the brick
may be any (coprime) triple of positive integers, the maximum number
of orientations of the placement of each further brick,
and also the maximum outdegree of each node of the state graph,
is six. (Consider three choices fixing the direction of the long edge
$t_n$ either forwards, sideways or up,
followed by another two choices for an orientation around the long axis.
If two edge lengths of the brick are equal, this set reduces to three choices.)
The encoding of the two-dimensional front of the brick insertions needs a matrix
of integers; we write one vector for each layer, separated by semicolons.

\subsection{$1\times 1\times 2$ Bricks}
\begin{figure}
\caption{
State diagram while tiling $2\times 2\times n$ rooms with
$1\times 1 \times 2$ dominoes.
}
\includegraphics[width=0.5\textwidth]{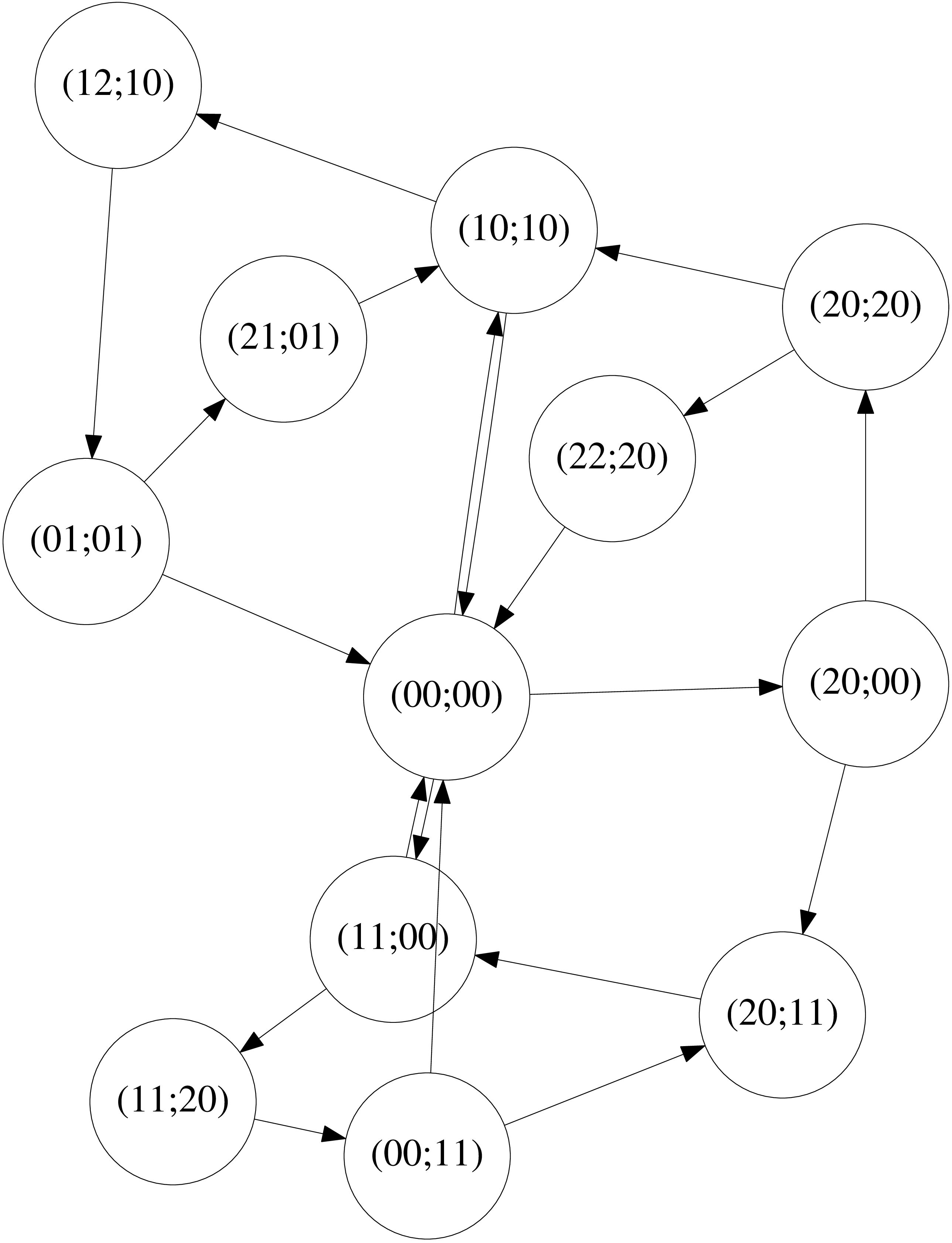}
\label{fig.state3d22211}
\end{figure}

To illustrate the notation, consider the simplest nontrivial
layout with rooms of $2\times 2$ cross sections filled with
bricks of size $1\times 1\times 2$. Figure \ref{fig.state3d22211}
evaluates to
\cite[A006253]{EIS}\cite{HockJMP24}
\begin{multline}
T_{1\times 1\times 2}(2,2,z)=
\frac { 1-z^2}
{ \left( 1+z^{2} \right)  \left( z^{4}-4z^{2}+1 \right) }
\\
=
1+2z^{2}+9z^{4}+32z^{6}+121z^{8}+450z^{10}+1681{
z}^{12}+6272z^{14}+23409z^{16}+87362z^{18}+\cdots;
\\
\hat T_{1\times 1\times 2}(2,2,z)=
\frac{z^2(2+3z^2-z^4)}{1-z^2}
=2z^{2}+5z^{4}+4z^{6}+4z^{8}+4z^{10}+4z^{12}+4
z^{14}+4z^{16}+\cdots.
\end{multline}
Filling $2\times 3 \times n$ rooms with $1\times 1\times 2$ bricks
generates a state diagram with 60 nodes and yields
\cite[A028447]{EIS}
\begin{multline}
T_{1\times 1\times 2}(2,3,z)=
-{\frac {14z^{12}-7z^{6}+1+z^{24}-7z^{18}+3z^{21}+16
z^{9}-3z^{3}-16z^{15}}{ \left( -7z^{3}+7z^{21}+48
z^{9}+1-13z^{6}+28z^{12}-48z^{15}-13z^{18}+z^{24}
 \right)  \left( z^{6}-z^{3}-1 \right) }}
\\
=
1+3z^{3}+32z^{6}+229z^{9}+1845z^{12}+14320z^{15}+
112485z^{18}+880163z^{21}+6895792z^{24}+54003765z^{27}+\cdots.
\end{multline}
Filling $3\times 3 \times n$ rooms with $1\times 1\times 2$ bricks
delivers
\cite[A028452]{EIS}
\begin{equation}
T_{1\times 1\times 2}(3,3,z)=
1+229z^{9}+117805z^{18}+64647289z^{27}+35669566217z^{
36}+\cdots.
\end{equation}
Filling $3\times 4 \times n$ rooms with $1\times 1\times 2$ bricks
generates a graph with 5544 nodes
\cite[A028453]{EIS},
and Lundow has also tabulated the
case of $4\times 4\times n$ rooms
\cite[A028454]{EIS}
\cite{Lundow1996}.

\subsection{$1\times 1\times 3$ Bricks}

\begin{figure}
\caption{
State diagram while tiling $2\times 3\times n$ rooms with
$1\times 1 \times 3$ 3-ominoes.
}
\includegraphics[width=0.6\textwidth]{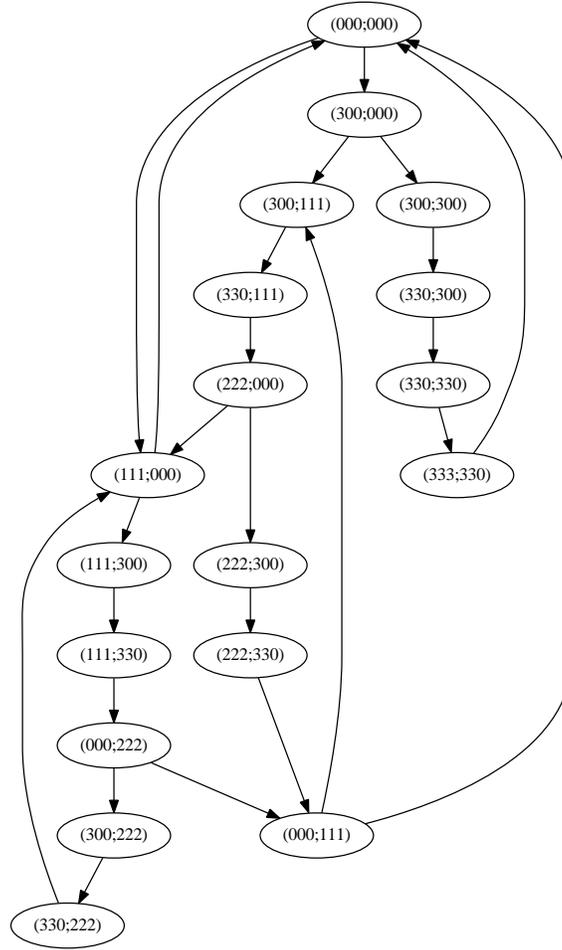}
\label{fig.state3d32311}
\end{figure}

Figure \ref{fig.state3d32311} describes placements of $1\times 1\times 3$ bricks
into $2\times 3\times n$ rooms
\cite[A233247]{EIS}:
\begin{equation}
T_{1\times 1\times 3}(2,3,z)=
{\frac {1-z^{6}-z^{4}}{ \left( z^{4}+1-z^{6} \right) 
 \left( 1-z^2-2z^4-z^6\right) }}
=
1+z^{2}+z^{4}+4z^{6}+9z^{8}+16z^{10}+36z^{12}+\cdots.
\label{eq.state3d32311}
\end{equation}
The expansion coefficients are $T_{1\times 1\times 3}(2,3,N)=T_{1 \times 3}^2(3,N)$,
the squares of those in \eqref{eq.state331}, because all solutions are classified
as stacks of two independent solutions of the 2-dimensional problem \cite{ReadAM24}.

Tiling $3\times 3\times n$ rooms with
$1\times 1 \times 3$ bricks yields
\cite[A233289]{EIS}
\begin{multline}
T_{1\times 1\times 3}(3,3,z)=
{\frac {z^{12}-z^{15}-4z^{9}-z^{21}-2z^{6}+z^{18}-z
^{3}+1}{ \left(-z^{27}-3z^{21}+z^{18}-6z^{15}-17z^{12}-
15z^{9}-2z^{6}-2z^{3}+1 \right)  \left(1-z^3 \right) 
 }}
=
\\
1+2z^{3}+4z^{6}+21z^{9}+92z^{12}+320z^{15}+1213{
z}^{18}+4822z^{21}+18556z^{24}+\cdots.
\end{multline}

Tiling $3\times 4\times n$ rooms with
$1\times 1 \times 3$ bricks yields \cite[A237355]{EIS}
\begin{equation}
T_{1\times 1\times 3}(3,4,z)
=
1
+3z^4
+9z^8
+92z^{12}
+749z^{16}
+4430z^{20}
+30076z^{24}
+217579z^{28}+\cdots.
\end{equation}
The full GF is not noted because the numerator is a polynomial of order 87 in $z^4$,
and the denominator a polynomial of order 88 in $z^4$.

\subsection{$1\times 2\times 2$ Bricks}

\begin{figure}
\caption{
State diagram while tiling $2\times 2\times n$ rooms with
$1\times 2 \times 2$ 4-ominoes.
}
\includegraphics[width=0.3\textwidth]{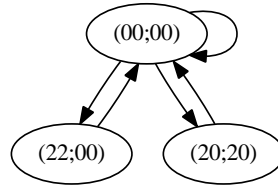}
\label{fig.state3d22221}
\end{figure}

Figure \ref{fig.state3d22221}
shows one loop of 1 step and two circuits of 2 steps \cite[A001045]{EIS}:
\begin{equation}
T_{1\times 2\times 2}(2,2,z)=
\frac{1}{1-z-2z^2}
=1+z+3z^{2}+5z^{3}+11z^{4}+21z^{5}+43z^{6}+85z
^{7}+171z^{8}+\cdots.
\end{equation}

\begin{figure}
\caption{
State diagram while tiling $2\times 3\times n$ rooms with
$1\times 2 \times 2$ 4-ominoes.
}
\includegraphics[width=0.7\textwidth]{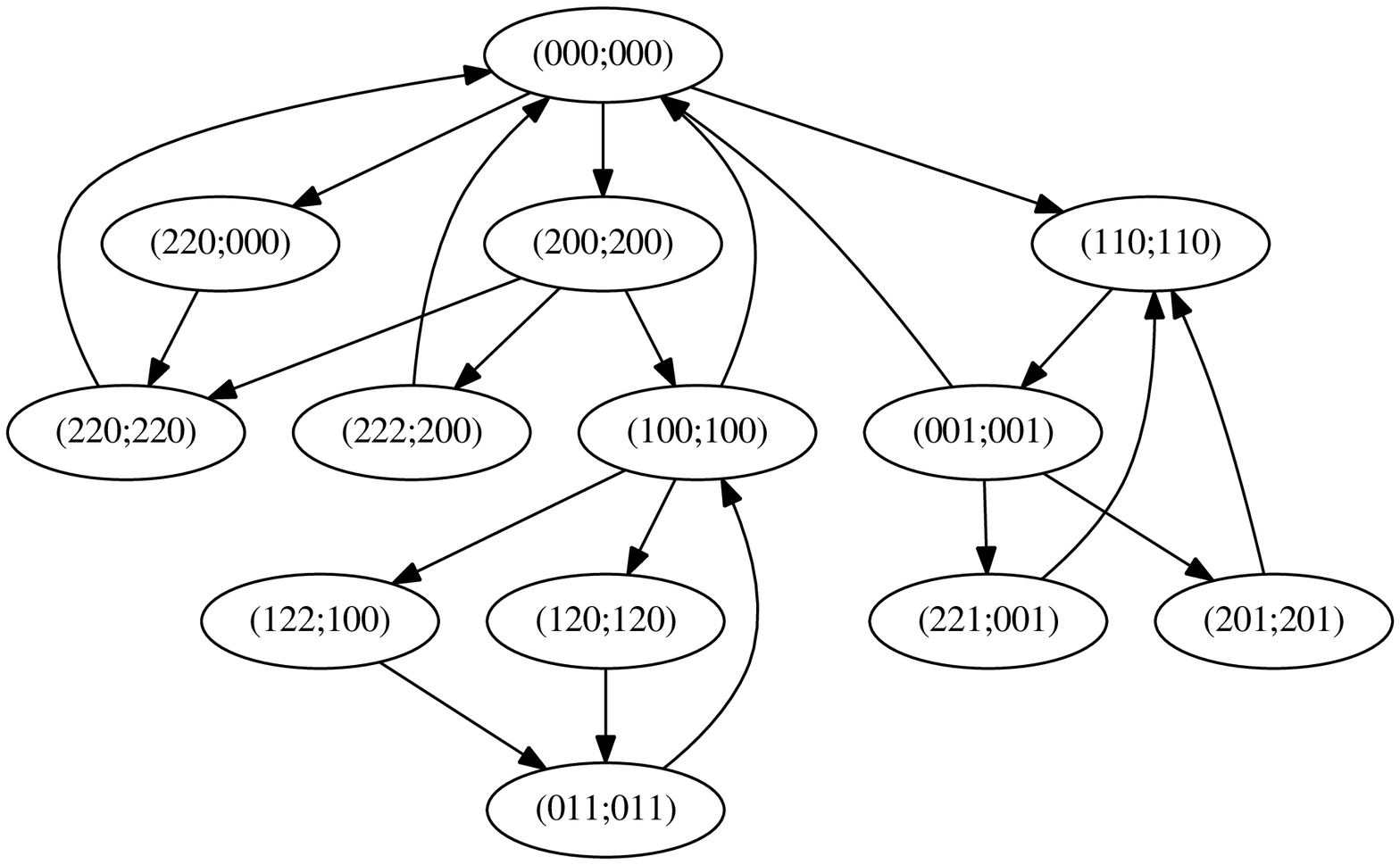}
\label{fig.state3d32221}
\end{figure}

Figure \ref{fig.state3d32221}
is \cite[A083066]{EIS}
\begin{equation}
T_{1\times 2\times 2}(2,3,z)=
{\frac {1-2z^{3}}{ \left(1- z^3 \right) \left(1- 6z^{3} \right) }}
=
1+5z^{3}+29z^{6}+173z^{9}+1037z^{12}+6221z^{15}+\cdots.
\end{equation}

\begin{figure}
\caption{
State diagram while tiling $3\times 3\times n$ rooms with
$1\times 2 \times 2$ 4-ominoes.
}
\includegraphics[width=0.95\textwidth]{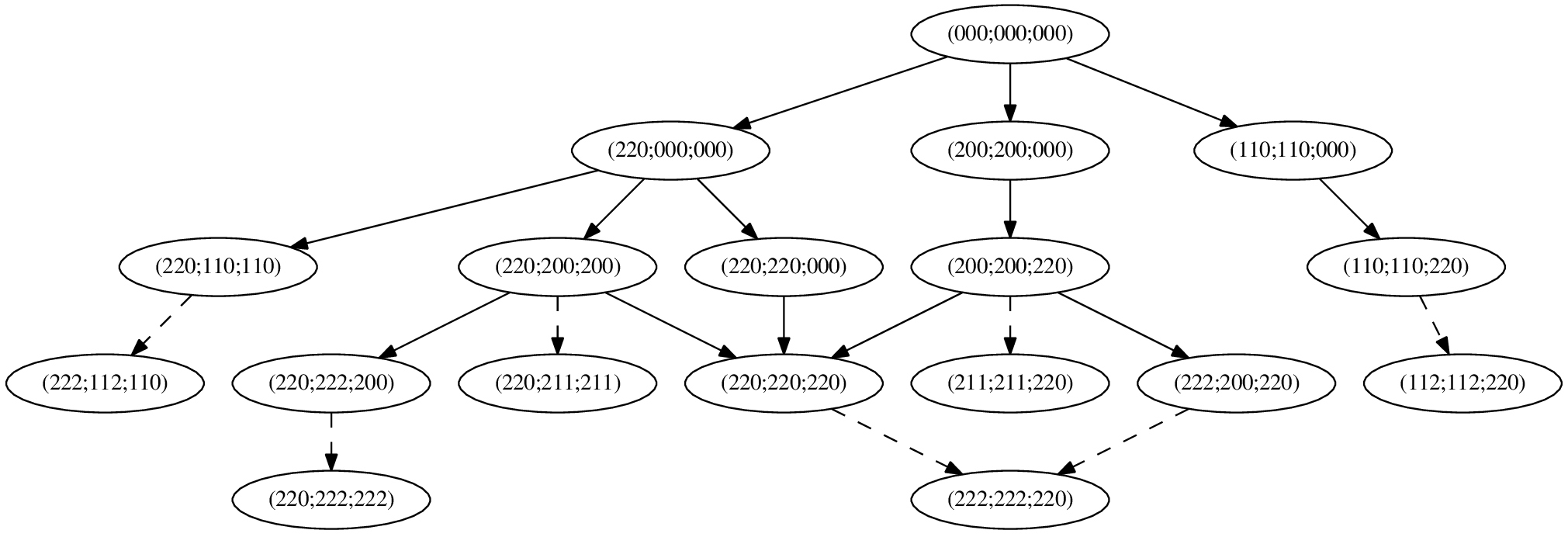}
\label{fig.state3d33221}
\end{figure}

Figure \ref{fig.state3d33221} is an example of a state diagram
without any of the required closed walks, $T_{1\times 2\times 2}(3,3,z)=1$.
This becomes plausible if one considers the subproblem of 
filling the first layer of cross section $3\times 3$ with any
combination of the $2\times 2$ or $1\times 2$ cross sections of
the brick---which cannot succeed because the area of 9 cannot
be partitioned into parts of 4 and 2\@.
This result is of a more general nature: if the cross section $m\times k$
of the room has an odd area $km$, at least one of the three cross sections
$t_mt_n$, $t_mt_k$ or $t_kt_n$ of the brick needs to be odd to allow
a complete filling.

Filling $3\times 4 \times n$ rooms with $1\times 2\times 2$
bricks is described by \cite[A237356]{EIS}
\begin{multline}
T_{1\times 2\times 2}(3,4,z)=
{\frac { \left( 1-2z^{6} \right)  \left(-120z^{18}+122z^{
12}-24z^{6}+1 \right) }{ \left(1- z^6 \right) 
\left( 2640{
z}^{24}-2540z^{18}+646z^{12}-54z^{6}+1 \right) }}
\\
=
1+29z^{6}+1065z^{12}+41097z^{18}+1602289z^{24}+
62603505z^{30}+2447085377z^{36}+\cdots.
\end{multline}

Filling $4\times 4 \times n$ rooms with $1\times 2\times 2$
bricks is described by
\begin{equation}
T_{1\times 2\times 2}(4,4,z)=
\frac{p_{1\times 2\times 2}(4,4,z)}
{q_{1\times 2\times 2}(4,4,z)}
=
1+z^{4}+165z^{8}+1065z^{12}+44913z^{16}+561689z^{20}
+\cdots,
\end{equation}
where
\begin{multline}
p_{1\times 2\times 2}(4,4,z)=
\left(1- 2z^{4} \right)  \left( 1+2z^{4} \right)
 ( 16896000z^{68}-21811200z^{64}-27278080z^{60}+
\\
43889536z^{56}+11612256z^{52}-32759456z^{48}+1863960z^{
44}+11174296z^{40}-2373860z^{36}
\\
-1742780z^{32}+538060z^
{28}+111540z^{24}-46358z^{20}-1940z^{16}+1612z^{12}-44
z^{8}-18z^{4}+1 );
\end{multline}
\begin{multline}
q_{1\times 2\times 2}(4,4,z)=
-207578z^{20}+3942932672z^{56
}-46651584z^{36}+1620z^{16}-597079616z^{60}
\\
+1-5117931520{
z}^{76}+229966968z^{44}-3208396800z^{80}-194z^{8}-418425808
z^{52}-19z^{4}
\\
+2433024000z^{84}+3948z^{12}+351702z^
{24}+4524176z^{28}+7975271424z^{72}-7861379200z^{64}\
\\
-
11992040z^{32}-1077496088z^{48}+3593770496z^{68}+160007540
z^{40}.
\end{multline}

\subsection{$1\times 1\times 4$ Bricks}

\begin{figure}
\caption{
State diagram while tiling $2\times 4\times n$ rooms with
$1\times 1 \times 4$ 4-ominoes.
}
\includegraphics[width=0.7\textwidth]{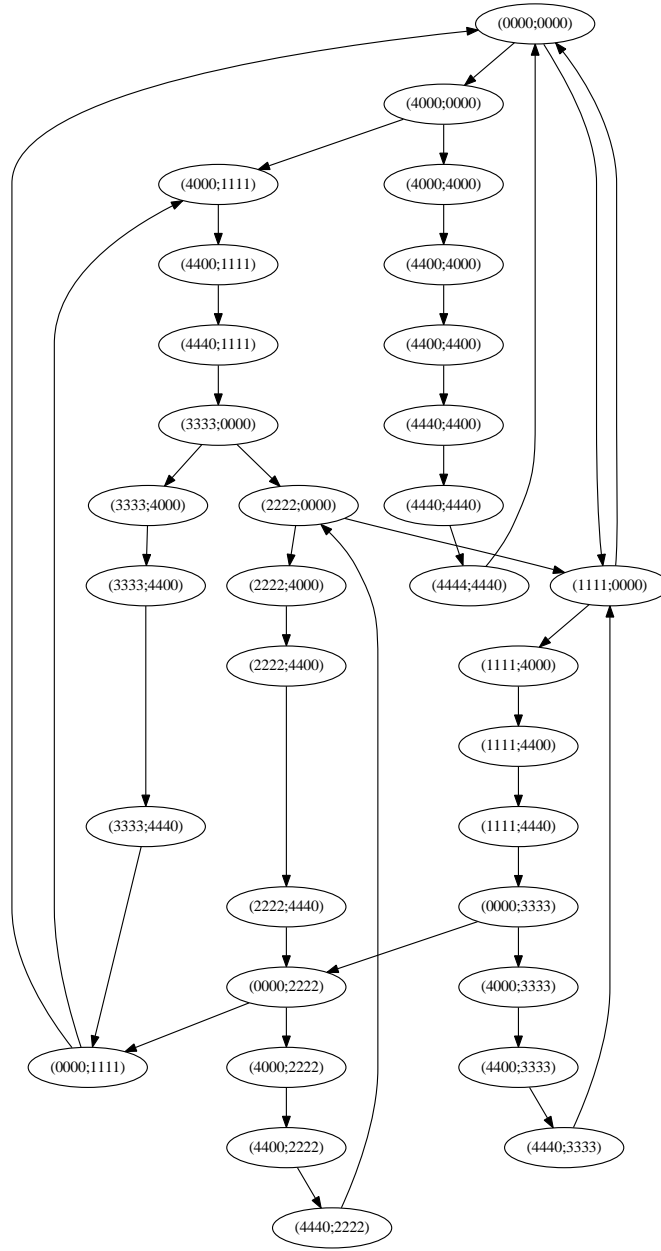}
\label{fig.state3d42411}
\end{figure}

Figure \ref{fig.state3d42411} constructs
\begin{multline}
T_{1\times 1\times 4}(2,4,z)=
{\frac {1-z^{4}-z^{8}+z^{12}-z^{6}}
{ \left(1-z^2-2z^4+z^8 \right) 
\left(-z^{12}+z^{4}+1-z^{8} +z^{6} \right) }}
\\
=
1+z^{2}+z^{4}+z^{6}+4z^{8}+9z^{10}+16z^{12}+25z^
{14}+49z^{16}+100z^{18}+196z^{20}+\cdots.
\end{multline}
The coefficients are the squares $T_{1\times 1 \times 4}(2,4,N)=T_{1\times 4}^2(4,N)$
of the coefficients of \eqref{eq.state314} for the reason
discussed in conjunction with Eq.\ \eqref{eq.state3d32311}.

Similarly the tiling $3\times 4\times n$ rooms with $1\times 1\times 4$ 4-ominoes
is counted by the cubes
$T_{1\times 1 \times 4}(3,4,N)=T_{1\times 4}^3(4,N)$:
\begin{equation}
T_{1\times 1\times 4}(3,4,z)=
\frac{p_{1\times 1\times 4}(3,4,z)}
{q_{1\times 1\times 4}(3,4,z)}
=
1+z^{3}+z^{6}+z^{9}+8z^{12}+27z^{15}+64z^{18}+125{
z}^{21}+343z^{24}+\cdots,
\end{equation}
where
\begin{multline}
p_{1\times 1\times 4}(3,4,z)=
z^{48}+z^{45}+2z^{39}-2z^{36}-3z^{33}-z^{30}
+5z^{27}+6z^{24}+3z^{21}
\\
+5z^{18}-6z^{15}-8z^{12
}-3z^{9}-2z^{6}+1;
\end{multline}
\begin{multline}
q_{1\times 1\times 4}(3,4,z)=
 \left(1-z^3-3z^6-3z^9-z^{12} \right)
 \left(-z^{12}-z^{9}+1 \right)
\\
\times
 \left(-z^{36}+3z^{33}-6z^{30}+7z^{27}-3z^{
24}-3z^{21}+4z^{18}-2z^{15}-4z^{12}+4z^{9}+z^{6}+
1 \right) .
\end{multline}

Tiling $4\times 4\times n$ rooms with $1\times 1\times 4$ 4-ominoes
is \cite[A233291]{EIS}
\begin{multline}
T_{1\times 1\times 4}(4,4,z)=
\frac{
p_{1\times 1\times 4}(4,4,z)}
{q_{1\times 1\times 4}(4,4,z)}
=
1+2z^{4}+4z^{8}+8z^{12}+45z^{16}+248z^{20}+1032{
z}^{24}+\cdots,
\end{multline}
with
\begin{multline}
p_{1\times 1\times 4}(4,4,z)=
1+30z^{28}-20z^{20}-112z^{80}+z^{116}-171z^
{40}-151z^{56}+90z^{32}
\\
+34z^{72}+69z^{36}+174z^{60}
+5z^{108}+61z^{24}+3z^{112}-4z^{12}-z^{4}+188z^{
64}-8z^{104}+63z^{76}
\\
-57z^{48}-166z^{68}-34z^{88}+
251z^{52}+z^{120}-3z^{100}-48z^{92}-104z^{44}-4z^
{8}-z^{124}+11z^{84}+39z^{96}-31z^{16} ;
\end{multline}
\begin{multline}
q_{1\times 1\times 4}(4,4,z)=
\left( 1+z^{4}
 \right)
( 1+184z^{28}+25z^{20}+89z^{80}+9z^{116}
+z^{136}-758z^{40}
\\
-715z^{56}+49z^{32}+436z^{72}+120
z^{36}+1435z^{60}-55z^{108}+72z^{24}+z^{112}+2z^{12
}-4z^{4}
\\
-1830z^{64}+120z^{104}-2z^{132}-883z^{76}-
1557z^{48}+545z^{68}-306z^{88}+917z^{52}+2z^{120}-
72z^{100}
\\
+z^{128}+75z^{92}+990z^{44}+2z^{8}-6z^{
124}+384z^{84}+120z^{96}-54z^{16}) .
\end{multline}

\subsection{$1\times 2\times 3$ Bricks}
Tiling $2\times 2\times n$ rooms with $1\times 2\times 3$ 6-ominoes
gives the powers of 2 \cite[A000079]{EIS},
\begin{multline}
T_{1\times 2\times 3}(2,2,z)=
\frac{1}{1-2z^2}
=
1+2z^{2}+4z^{4}+8z^{6}+16z^{8}+32z^{10}+\cdots.
\end{multline}
Results of that type are generally understood by the constraint that
the $2\times 2$ cross section of the room permits only two different orientations
of the $1\times 2$ cross section of the brick, and that after each
such placement the next placement is forced to return to the straight
profile $(00\ldots0)$.

\begin{figure}
\caption{
State diagram while tiling $2\times 3\times n$ rooms with
$1\times 2 \times 3$ 6-ominoes.
}
\includegraphics[width=0.6\textwidth]{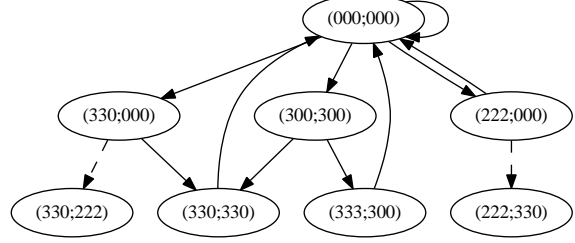}
\label{fig.state3d32123}
\end{figure}

Tiling $2\times 3\times n$ rooms with $1\times 2\times 3$ 6-ominoes
is analyzed in Figure \ref{fig.state3d32123} \cite[A103143]{EIS}:
\begin{multline}
T_{1\times 2\times 3}(2,3,z)=
\frac{1}{1-z-z^2-3z^3}
=
1+z+2z^{2}+6z^{3}+11z^{4}+23z^{5}+52z^{6}+108z
^{7}+\cdots.
\end{multline}
After pruning $(330;222)$ and $(222;330)$,
the single loop, the circuit with 2 steps passing through (222;000)
and the 3 circuits with 3 steps that contribute to the denominator of
this generating function are easily recognized in the figure.

Tiling $3\times 3\times n$ rooms with $1\times 2\times 3$ 6-ominoes
is counted by \cite[A237357]{EIS}
\begin{multline}
T_{1\times 2\times 3}(3,3,z)=
\frac { 1-z^3 }{1-22z^{6}- 7z^{3}-36z^{9}}
=
1+6z^{3}+64z^{6}+616z^{9}+5936z^{12}+57408z^{15}+\cdots;
\\
\hat T_{1\times 2\times 3}(3,3,z)=
\frac {2z^{3} \left( 3+11z^{3}+18z^{6} \right) }
{ 1-z^3}
=
6z^{3}+28z^{6}+64z^{9}+64z^{12}+64z^{15}+64z^{
18}+\cdots.
\end{multline}

Tiling $3\times 4\times n$ rooms with $1\times 2\times 3$ 6-ominoes
is summarized by the GF \cite[A237358]{EIS}
\begin{equation}
T_{1\times 2\times 3}(3,4,z)=\frac{p_{1\times 2\times 3}(3,4,z)}{q_{1\times 2\times 3}(3,4,z)}
=
1+z^{2}+11z^{4}+64z^{6}+296z^{8}+1716z^{10}+9123{z
}^{12}+\cdots
\end{equation}
with numerator
\begin{multline}
p_{1\times 2\times 3}(3,4,z)=
{\left(1-z^2\right)  \left( 1+z^{2}
 \right)  \left(1- 3z^{2} \right)  \left( 3z^{4}+2z^{2}+1
 \right)  \left( 1-z^{4}-7z^{6}+9z^{12} \right) }
\end{multline}
and denominator
\begin{multline}
q_{1\times 2\times 3}(3,4,z)=
504z^{ 12}+306z^{22}+1-1012z^{18}+103z^{14}-2z^{2}+54z^{32
}-162z^{34}
\\
-450z^{28}+74z^{24}-14z^{4}-487z^{16}-42
z^{6}-448z^{20}+915z^{26}+237z^{10}+873z^{30}+42{
z}^{8}.
\end{multline}

\section{Summary}
We have transformed Read's profile vectors of incomplete tilings
into Transfer Matrices of associated digraphs, and obtained generating functions
for some tilings of rectangular floors and rooms by 2- and 3-dimensional rectangular
tiles.

\bibliographystyle{amsplain}
\bibliography{all}

\end{document}